\documentclass[12pt]{amsart}
\usepackage{amscd,amssymb,epsfig}
\calclayout

\numberwithin{equation}{section}

\newcommand\cof{\operatorname{cof}}

\newcommand{\ud}{d}

\renewcommand{\div}{\operatorname{div}}

\newtheorem{thm}{Theorem}[section] 
 
\newtheorem{lemma}[thm]{Lemma}

\newtheorem{rem}[thm]{Remark}

\def\whsq{\vbox to 5.8pt 
{\offinterlineskip\hrule 
\hbox to 5.8pt{\vrule height 
5.1pt\hss\vrule height 5.1pt}\hrule}}

\def\<{\langle} 
\def\>{\rangle} 
\def\PP{{\mathop{{\rm I}\kern-.2em{\rm P}}\nolimits}} 
\def\FF{{\mathop{{\rm I}\kern-.2em{\rm F}}\nolimits}}   
\def\ZZ{{\mathop{{\rm I}\kern-.2em{\rm Z}}\nolimits}} 
 
\newlength{\sidemargin} 
\begin{document}
\title[]{
Spline element method for Monge-Amp\`ere equations
}

\thanks{The author was supported in part by NSF grant DMS-0811052 and the Sloan Foundation. }
\author{Gerard Awanou}
\address{Department of Mathematics, Statistics, and Computer Science, M/C 249.
University of Illinois at Chicago, 
Chicago, IL 60607-7045, USA}
\email{awanou@uic.edu}  
\urladdr{http://www.math.uic.edu/\~{}awanou}

\maketitle

\begin{abstract}
We analyze the convergence of an iterative method for solving the nonlinear system resulting from a natural discretization of the Monge-Amp\`ere equation with $C^1$ conforming approximations. We make the assumption, supported by numerical experiments for the two dimensional problem, that the discrete problem has a convex solution. The method we analyze is the discrete version of Newton's method in the vanishing moment methodology. 
Numerical experiments are given in the framework of the spline element method. 
\end{abstract}

\section{Introduction}
This paper addresses the numerical solution of the Dirichlet problem for the Monge-Amp\`ere 
equation
\begin{equation}
\det  D^2 u = f  \, \text{in} \, \Omega, \quad u=g \, \text{on} \, \partial \Omega. 
\label{m1}
\end{equation}
Here $D^2 u=\bigg( \frac{\partial^2 u}{\partial x_i \partial x_j}\bigg)_{i,j=1,\ldots, n}  $ 
is the Hessian of $u$
and $f,g$ are given functions with $f\geq c_0 > 0$. The domain $\Omega \subset \mathbb{R}^n, n=2,3$ is assumed to be bounded and convex with a polygonal boundary and $\partial \Omega$
denotes its boundary. 

Let $V_h$ denote a finite dimensional space of $C^1$ functions which are piecewise polynomials of degree $d \geq 2$ and let us assume that $f \in L^1(\Omega)$. We consider the discrete problem: find $u_h \in V_h$ such that
\begin{align}
\begin{split}
\int_{\Omega} v_h \det D^2 u_h \ud x & = \int_{\Omega} f v_h \ud x, \forall v_h \in V_h \cap H_0^1(\Omega) \\
u_h & = g_h \text{on} \, \partial \Omega,
\end{split} \label{m1h}
\end{align}
where $g_h$ is the natural interpolant in $V_h$ of a smooth extension of $g$. In this paper, we make the assumption that \eqref{m1h} has a strict convex solution $u_h$. We analyze the convergence of 
the following iterative method. Given an initial guess $u^0_h \in V_h$ with $u_h^{0}  = g_h$ $ \text{on} \, \partial \Omega$, find $u^{k+1}_h \in V_h$ such that $u_h^{k+1}  = g_h$ $\text{on} \, \partial \Omega$ and for $\epsilon >0$
\begin{align}\label{newtonvm}
\begin{split}
\epsilon& \int_{\Omega} \Delta u_h^{k+1} \Delta v_h \ud x  +  \int_{\Omega} [(\cof D^2 u_h^k) D u_h^{k+1}] \cdot D v_h  \ud x  =  -\int_{\Omega} f v_h  \ud x \\
& \quad \quad + \epsilon^3 \int_{\partial \Omega} \frac{\partial v_h}{\partial n} \ud s  
+   \frac{n-1}{n}  \int_{\Omega} [(\cof D^2 u_h^k) D u_h^{k}] \cdot D v_h  \ud x,  \forall v_h \in V_h \cap H_0^1(\Omega).
\end{split}
\end{align}

The main difficulties of the numerical resolution of \eqref{m1} is that when it does not have a smooth solution, Newton's method (i.e. \eqref{newtonvm} with $\epsilon=0$) breaks down. 

 In \cite{AwanouPseudo10} we show that \eqref{m1h} is well defined and has a strict convex solution when \eqref{m1} has a smooth strictly convex solution. Less restrictive conditions under which \eqref{m1h} has a strict convex solution are addressed in \cite{Awanou-Std04} in the framework of the Aleksandrov theory of the Monge-Amp\`ere equation. The assumption of existence of a strictly convex solution of \eqref{m1h} is supported in this paper by numerical experiments in two dimension. 
 We prove the convergence of the iterations \eqref{newtonvm} to a limit $u_{\epsilon,h}$ which solves a discrete variational problem. With that result, one may prove a quadratic convergence rate for \eqref{newtonvm} as an iterative method converging to $u_{\epsilon,h}$, using for example the techniques of \cite{AwanouPseudo10}. That issue is not addressed in this paper since \eqref{newtonvm} is not a direct method for solving \eqref{m1h}.

For $C^1$ conforming approximations, we use the spline element method \cite{Awanou2003, Awanou2005a, Awanou2006, Baramidze2006, Hu2007,Awanou2008}. It
uses piecewise polynomials of arbitrary degree
and Lagrange multipliers to enforce continuity and smoothness conditions as well as constraints.
However, unlike other methods which also use Lagrange multipliers, the constraints here
are enforced exactly. More details are given in Section \ref{splineS}. An alternative to the spline element method is the Argyris finite element for the two dimensional problem or concepts from isogeometric analysis \cite{Vuong10}. The study of $C^1$ conforming approximations provides a natural setting for presenting techniques for proving results on the numerical analysis of Monge-Amp\`ere equations. These techniques may be extended to the setting of isogeometric analysis, mixed finite elements, Lagrange elements or the standard finite difference method.

The iterative method \eqref{newtonvm} is the discrete version of Newton's method in the vanishing moment methodology. In \cite{Feng2009a,Feng2009b,Feng2009,Neilan2010}, it was proposed to solve \eqref{m1} by the singular perturbation problem
\begin{align*} 
- \epsilon \Delta^2 u + \text{det} \ D^2 u = f,
\ \text{in} \ \Omega, \quad u=g, \ \Delta u = \epsilon^2 \ \text{on} \ \partial \Omega.
\end{align*}
The point of view we take in this paper is rather different. We do not view \eqref{newtonvm} as the discretization of a partial differential equation at the continuous level. But as an iterative method for solving the discrete nonlinear system of equations \eqref{m1h}. The results we prove may be viewed as partial discrete analogues of the assumptions made in \cite{Feng2009a,Feng2009b,Feng2009,Neilan2010} and partially proved in \cite{Feng2012}. It should be noted that our approach can also be reproduced at the continuous level. Combined with the Aleksandrov theory of the Monge-Amp\`ere equation, continuous analogues of the results presented in this paper may be proven. We wish to address  in a separate work this new approach to the vanishing moment methodology.

Existing numerical work on Monge-Amp\`ere type equations can be divided into three main categories. In the first category, the
Monge-Amp\`ere equation is treated as a nonlinear elliptic equation and the impressive tools of classical finite element analysis
are used to derive error estimates and convergence of iterative methods at the continuous and discrete levels
\cite{Loeper2005,Belgacem2006,Bohmer2008,Feng2009,Feng2009a,Bohmer2010,Brenner2010b,Brenner2010a,AwanouPseudo10}. In the second category we have methods which are illustrated to 
approximate the weak solutions (some of them considering only smooth solutions) although there is no evidence of theoretical convergence \cite{Dean2003,Dean2004,Dean2006,Mohammadi2007,GlowinskiICIAM07,Feng2009b,Benamou2010,Oberman2010b,Zheligovsky10,AwanouPseudo10,Lakkis11b}. In the final category, we have 
methods which are proven to converge to weak solutions. Here, since the solution may not be regular, one cannot expect to have a lot of information
about convergence rates \cite{Oliker1988, Kochengin1998, Caffarelli1999, Wang2004, Oberman2008,Oberman2010a}.
This paper falls mainly in the first two categories. 

The paper is organized as follows: in the second section, we introduce some notation and give some preliminary results. In section \ref{vanishing} we prove the convergence of Newton's method in the vanishing moment methodology.  
The last section is devoted to numerical experiments.

\section{Notation and Preliminaries}
We use the usual notation $L^p(\Omega), 1 \leq p \leq \infty$ for the Lebesgue spaces and $W^{k,p}(\Omega)$ for the Sobolev spaces with norms $||.||_{k,p}$ and semi-norm $|.|_{k,p}$. In particular,
$H^k(\Omega)=W^{k,2}(\Omega)$ and in this case, the norm and semi-norms will be denoted respectively by 
$||.||_{k}$ and semi-norm $|.|_{k}$. For two $n \times n$ matrices  $A, B$, we recall the Frobenius inner product $A:B = \sum_{i,j=1}^n A_{ij} B_{ij}$, where $A_{ij}$ and $B_{ij} $ refer to the entries of  the corresponding matrices. For a matrix field $A$, we denote by $\div A$ the vector obtained by taking the divergence of each row. We will use the notation
$$
||A||_{\infty} : = \max_{i,j} |a_{ij}|,
$$
for a matrix $A=(a_{ij})_{i,j=1,\ldots,n}$ 
and denote by $n$ the unit outward normal vector to $\partial \Omega$. 

We require our approximation spaces $V_h$ to satisfy the following properties:
There exists an interpolation operator
$I_h$ mapping $W^{l+1,p}(\Omega)$ into the space $V_h$ for $1 \leq p \leq \infty, 0 \leq l \leq d$ 
such that 
\begin{equation}
|| v -I_h v ||_{k,p} \leq C_{} h^{l+1-k} ||v||_{l+1,p}, \label{schum}
\end{equation}
for $0 \leq k \leq l$ 
and the inverse estimates
\begin{equation}
||v||_{s,p} \leq C_{} h^{l-s+\text{min}(0,\frac{n}{p}-\frac{n}{q})} ||v||_{l,q}, \forall v \in V_h, \label{inverse}
\end{equation}
and for $0 \leq l \leq s, 1 \leq p,q\leq \infty$.

The above assumptions are known to be satisfied for standard finite element spaces \cite{Brenner02}. For the spline spaces used in the computations, \eqref{schum} is known to hold \cite{Lai2007}. One may view \eqref{inverse} as a consequence of Markov inequality, \cite{Lai2007} p. 2. One can then prove global estimates 
as in \cite{Lai98}.

It follows from \eqref{schum} that
\begin{equation} \label{stable}
|| I_h v ||_{k,p} \leq C ||  v ||_{k,p}, 
\end{equation}
for $1 \leq p \leq \infty$ and $0 \leq k \leq d$.

We will need the following lemma whose proof can be found in  \cite{AwanouPseudo10}
\begin{lemma} \label{det-lem} We have
\begin{equation}\label{detform}
\det D^2 v = \frac{1}{n} (\cof D^2 v): D^2 v =
\frac{1}{n} \div\big((\cof D^2 v) D v \big). 
\end{equation}
And for $F(v) = \det D^2 v $ we have
$$
F'(v) (w) =  (\cof D^2 v): D^2 w = \div \big((\cof D^2 v) D w \big), 
$$
 for $v,w$ 
sufficiently smooth.
\end{lemma}
We make the usual convention of denoting  constants by $C$. Let us denote by $\lambda_1(D^2 v)$ and $\lambda_n(D^2 v)$ the smallest and largest eigenvalues respectively of $D^2 v$. We will use the notation $D v$ to denote the gradient vector of the function $v$ and recall that $\cof A$ denotes the matrix of cofactors of the matrix $A$.

We make in this paper the assumption that \eqref{m1h} has a strictly convex solution  $u_h$ with $0< 2C_0 \leq \lambda_1(D^2 u_h) \leq \lambda_n(D^2 u_h) \leq C_4$ for constants $C_0$ and $C_4$ independent of $h$. We define for $\rho >0$
$$
B_{\rho}(u_h) = \{ \, v_h \in V_h, ||v_h-u_h||_1 \leq \rho\, \}.
$$
By the continuity of the eigenvalues of a matrix as a function of its entries, $v_h \in B_{\rho}(u_h)$ is strictly convex for $\rho$ sufficiently small. Moreover, $\lambda_1(D^2 v_h) \geq C_0$ for $v_h \in B_{\rho}(u_h)$. See \cite{AwanouPseudo10} for details.

\section{Convergence of the discrete vanishing moment methodology} \label{vanishing}

By Lemma \ref{det-lem}, we have for $w_h \in B_{\rho}(u_h)$ and $ v_h \in V_h \cap H_0^1(\Omega)$,
\begin{align} \label{cof-det-lem}
\begin{split}
 \int_{\Omega} [(\cof D^2 w_h ) D w_h ] \cdot D v_h  \ud x &= -  \int_{\Omega} \div [(\cof D^2 w_h ) D w_h] v_h  \ud x \\
 & = -n  \int_{\Omega}  (\det D^2 w_h) v_h  \ud x.
 \end{split}
\end{align}
Thus, we can rewrite \eqref{newtonvm} as
\begin{align}\label{newtonvm1}
\begin{split}
\epsilon \int_{\Omega} \Delta u_h^{k+1} \Delta v_h \ud x  +  \int_{\Omega} [(\cof D^2 u_h^k) D u_h^{k+1}] \cdot D v_h  \ud x   =  \epsilon^3 \int_{\partial \Omega} \frac{\partial v_h}{\partial n} \ud s   
+\int_{\Omega} p_h^k v_h  \ud x,
\end{split}
\end{align}
 for all $v_h \in V_h \cap H_0^1(\Omega)$ with
\begin{equation} \label{p-form}
 p_h^k=-f -(n-1) \det  D^2 u_h^k.
 \end{equation}

Given $u_h^k \in B_{\rho}(u_h)$, with $u_h^k=g_h$ on $\partial \Omega$, let $\hat{u}_h^{k+1}$ satisfy $\hat{u}_h^{k+1}=g_h$ on $\partial \Omega$ and for all $ v_h \in V_h \cap H_0^1(\Omega)$,
\begin{equation} \label{newton}
  \int_{\Omega} [(\cof D^2 u_h^k) D \hat{u}_h^{k+1}] \cdot D v_h  \ud x  =  \int_{\Omega} p_h^k v_h  \ud x. 
\end{equation}
We claim that there exists a constant $C_1 < 1$ such that
\begin{equation} \label{newton-step}
|| \hat{u}_h^{k+1} - u_h ||_1 \leq C_1 ||u_h^k-u_h||_1.
\end{equation}
We also claim that for $\epsilon$ sufficiently small
\begin{equation} \label{vm-step}
 || \hat{u}_h^{k+1} - u_h^{k+1} ||_1 \leq (1 - C_1) \rho.
\end{equation}
The proof of the above claims are given below. Let us first derive some consequences. We would then have
$$
 || u_h^{k+1} -u_h ||_1 \leq || \hat{u}_h^{k+1} - u_h^{k+1} ||_1 +   || \hat{u}_h^{k+1} - u_h ||_1 \leq  C_1 \rho +  \rho - C_1 \rho \leq \rho. 
$$
We conclude that given an initial guess $u_h^0$ in $B_{\rho}(u_h)$, we have $ u_h^k \in B_{\rho}( u_h)$ for all $k$. 
Therefore, there exists a subsequence, which is also denoted $u_h^k$, which converges to an element 
$u_{\epsilon,h}$ of $ B_{\rho}(u_h)$.  We claim that the convex function $u_{\epsilon,h} \in V_h$ solves $u_{\epsilon,h}  = g_h$ $\text{on} \, \partial \Omega$ and for all $ v_h \in V_h \cap H_0^1(\Omega)$,
\begin{align}\label{vm}
\begin{split}
\epsilon& \int_{\Omega} \Delta u_{\epsilon,h} \Delta v_h \ud x  -  \int_{\Omega} (\det D^2 u_{\epsilon,h}) v_h  \ud x  =  -\int_{\Omega} f v_h  \ud x + \epsilon^3 \int_{\partial \Omega} \frac{\partial v_h}{\partial n} \ud s.  
\end{split}
\end{align}
Since $\rho > 0$ is taken small but otherwise arbitrary, we conclude that $|| u_{\epsilon,h}-u_h||_1 \to 0$ as $h \to 0$.

\subsection{Proof of \eqref{newton-step}}
The equation \eqref{newton-step} is nothing but a Newton's step for $C^1$ conforming approximations of the Monge-Amp\`ere equation. 
Starting with $u_h^k$, one step of Newton's method applied to \eqref{m1h} produces $ \hat{u}_h^{k+1}$ as a solution of
\begin{align*}
\begin{split}
\int_{\Omega} [\div ( (\cof D^2 u_h^k) D (\hat{u}_h^{k+1} -u_h^k) ] v_h  \ud x & = - \int_{\Omega} (\det D^2 u_h^k - f) v_h \ud x, \forall v_h \in V_h \cap H_0^1(\Omega) \\
\hat{u}_h^{k+1} & = g_h \text{on} \, \partial \Omega,
\end{split} 
\end{align*}
Using integration by parts and taking into account \eqref{cof-det-lem}, we obtain \eqref{newton-step}.

From \cite{AwanouPseudo10}, we know that for $h$ sufficiently small, there exists a constant $C <1$ such that
$$
|| \hat{u}_h^{k+1} - u_h ||_1 \leq C ||u_h^k-u_h||_1^2 \leq C \rho ||u_h^k-u_h||_1.  
$$
Thus with $\rho \leq 1$ and $C_0 = C \rho$, we obtain the result.

\subsection{Proof of \eqref{vm-step}} We view this step as a correction with the regularization. 

Substituting \eqref{newton} into \eqref{newtonvm1}, we obtain
\begin{multline*}
\epsilon \int_{\Omega} (\Delta u_h^{k+1}-\Delta \hat{u}_h^{k+1})  \Delta v_h \ud x  +  \int_{\Omega} [(\cof D^2 u_h^k) D (u_h^{k+1}-  \hat{u}_h^{k+1} )] \cdot D v_h  \ud x  \\ =  \epsilon^3 \int_{\partial \Omega} \frac{\partial v_h}{\partial n} \ud s  - \epsilon \int_{\Omega} \Delta \hat{u}_h^{k+1} \Delta v_h \ud x.
\end{multline*}
We now use various letters $C_i, i= 2, 3, 4$ to update various constants.

Substituting $v_h=u_h^{k+1}-  \hat{u}_h^{k+1}$ in the above equation, and using the strict convexity of $u^k_h$, we obtain using a trace estimate and inverse inequalities
\begin{align*}
\epsilon ||\Delta(u_h^{k+1}-  \hat{u}_h^{k+1})||_0^2 + C_2  |u_h^{k+1}-  \hat{u}_h^{k+1}|_1^2 & \leq C_3 \epsilon^3 \bigg( \int_{\partial \Omega}\bigg| \frac{\partial v_h}{\partial n} \bigg|^2 \ud s \bigg)^{\frac{1}{2}} \\
& \qquad \qquad  + \epsilon ||\Delta \hat{u}_h^{k+1} ||_0  ||\Delta v_h ||_0 \\
& \leq  C_3 \epsilon^3 ||v_h||_2  + C_4\epsilon ||\Delta \hat{u}_h^{k+1} ||_0  ||\Delta v_h ||_0 \\
& \leq  C_3 h^{-1} \epsilon^3 ||v_h||_1  + C_4 \epsilon h^{-2} || \hat{u}_h^{k+1} ||_1  || v_h ||_1.
\end{align*}
We conclude that
\begin{align*}
||u_h^{k+1}-  \hat{u}_h^{k+1}||_1 \leq C_3 h^{-1} \epsilon^3 + C_4 \epsilon h^{-2} (\rho + ||u_h||_1).
\end{align*}
We obtain \eqref{vm-step} if we choose $\epsilon$ such that
\begin{align*}
\epsilon \leq \min \bigg\{ \, \bigg(\frac{(1-C_1) h \rho }{3 C_3 } \bigg)^{\frac{1}{3}}, \frac{(1-C_1) h^2}{3 C_4}, \frac{(1-C_1) h^2\rho }{3 C_4 ||u_h||_1} \, \bigg\}.
\end{align*}

\subsection{Passage to the limit}
By an inverse estimate or the equivalence of norms in a finite dimensional space, the sequence $u_h^k$ is also bounded in $W^{2,n}(\Omega)$ and hence converges (up to a subsequence) in $W^{2,n}(\Omega)$ to a limit $u_{\epsilon,h}$. Passing in the limit in \eqref{newtonvm}, we obtain \eqref{vm} as follows. For $ v_h \in V_h \cap H_0^1(\Omega)$, we have
\begin{align*}
\bigg| \int_{\Omega}( \Delta u_{\epsilon,h} - \Delta u_h^{k+1})\Delta v_h \ud x \bigg| & \leq ||  \Delta u_{\epsilon,h} - \Delta u_h^{k+1} ||_0 ||\Delta v_h||_0 \leq C || u_{\epsilon,h} - u_h^{k+1}||_2 ||v_h||_2\\
& \to 0 \, \text{as} \, k \to \infty.
\end{align*}
Put 
$$
A_1 = \int_{\Omega} [(\cof D^2 u_h^k - \cof D^2 u_{\epsilon,h}) D u_h^{k+1}] \cdot D v_h  \ud x,
$$
and 
$$
A_2 =  \int_{\Omega} [ (\cof D^2 u_{\epsilon,h}) (D u_h^{k+1} - D u_{\epsilon,h} ) ] \cdot D v_h  \ud x.
$$
We have by Cauchy-Schwarz inequality and the inverse estimate \eqref{inverse}
\begin{align*}
|A_2|&  \leq C ||u_{\epsilon,h}||_{2,\infty}^{n-1} || u_h^{k+1} -  u_{\epsilon,h} ||_1 ||v_h||_1\\
&  \leq C h^{-(n-1)(2+\frac{n}{2})}  ||u_{\epsilon,h}||_2 || u_h^{k+1} -  u_{\epsilon,h} ||_1 ||v_h||_1\\
& \to 0 \, \text{as} \, k \to \infty.
\end{align*}
Let us denote by $(\cof)'$ the Fr\'echet derivative of the mapping $A \to \cof A$. Since $(\cof)'(A)(B)$ is the sum of terms which are products of $n-2$ components of $A$ and is linear in the components of $B$, we have
$$
||(\cof)'(D^2 v)(D^2 w)||_{0,\infty} \leq C ||D^2 v||_{2,\infty}^{n-2} ||D^2 w||_{2,\infty}.
$$
It follows that
\begin{align*}
|A_1| & \leq C \sum_{K \in \mathcal{T}_h} || u_h^k -  u_{\epsilon,h}||_{2,\infty} ||u_h^{k+1}||_{1,K} ||v_h||_{1,K} \\
& \leq C || u_h^k -  u_{\epsilon,h}||_{2,\infty}  ||u_h^{k+1}||_{1} ||v_h||_{1} \\
& \leq C h^{-(2+\frac{n}{2})} || u_h^k -  u_{\epsilon,h}||_2 ||u_h^{k+1}||_{1} ||v_h||_{1}\\
& \to 0 \, \text{as} \, k \to \infty,
\end{align*}
since the convergent sequence $||u_h^{k+1}||_{1}$ is bounded.

Finally
\begin{align*}
\bigg|  \int_{\Omega} [(\cof D^2 u_h^k) D u_h^{k+1}] \cdot D v_h  \ud x -  \int_{\Omega} [(\cof D^2 u_{\epsilon,h}) D u_{\epsilon,h}] \cdot D v_h  \ud x \bigg| & = |A_1+A_2| \\
&  \to 0 \, \text{as} \, k \to \infty.
\end{align*}
Passing in the limit in \eqref{newtonvm}, we have
\begin{align*} 
\begin{split}
\epsilon& \int_{\Omega} \Delta u_{\epsilon,h} \Delta v_h \ud x  +  \int_{\Omega} [(\cof D^2 u_{\epsilon,h}) D u_{\epsilon,h} ] \cdot D v_h  \ud x  =  -\int_{\Omega} f v_h  \ud x \\
& \quad \quad + \epsilon^3 \int_{\partial \Omega} \frac{\partial v_h}{\partial n} \ud s  
+   \frac{n-1}{n}  \int_{\Omega} [(\cof D^2 u_{\epsilon,h}) D u_{\epsilon,h}] \cdot D v_h  \ud x,  
\end{split}
\end{align*}
By \eqref{cof-det-lem} we obtain \eqref{vm}.

\subsection{Pointwise convergence of boundary data}
Since $u_h^{k+1}=g_h$ on $\partial \Omega$, it follows that $u_h^{k+1}$ is bounded on $\partial \Omega$. Passing to a subsequence, we conclude that $u_{\epsilon,h}  = g_h$ $\text{on} \, \partial \Omega$ as well.

We have proved the following theorem
\begin{thm}
Up to a subsequence the sequence defined by \eqref{newtonvm} converges to the solution $u_h$ of \eqref{m1h} for $\epsilon$ and $h$ sufficiently small and a sufficiently close initial guess. Moreover, as $k \to \infty$, the subsequence  converges to a convex solution of \eqref{vm}.
\end{thm}
\begin{rem}
Since $u_{\epsilon,h} \in B_{\rho}(u_h)$, $u_{\epsilon,h}$ is also convex. The convexity of the solution of the discrete variational problem obtained in the vanishing moment methodology, namely \eqref{vm}, has long been an open problem.
\end{rem}
We now prove that the whole sequence defined by \eqref{newtonvm} converges to a solution of  \eqref{vm}. For this, it is enough to prove that  \eqref{vm} has a unique solution in $B_{\rho}(u_h)$. We recall from \cite{AwanouPseudo10} that for $h$ sufficiently small and all $v_h \in B_{\rho}(u_h)$
\begin{equation} \label{p-d-K}
m |w|_{1,K}^2 \leq \int_{K} [(\cof\, D^2 v_h(x)) D w(x)] \cdot D w(x) \, \ud x \leq  M |w|_{1,K}^2, w \in H^1(K), 
\end{equation}
with constants $m$ and $M$ independent of $h$. To prove \eqref{p-d-K} one uses the positive definiteness of $\cof\, D^2 v_h$.

Let $u_{\epsilon,h}, v_{\epsilon,h}$ be two solution of \eqref{vm} in $B_{\rho}(u_h)$. For all $t \in [0,1]$, $t u_{\epsilon,h} + (1-t) v_{\epsilon,h} \in B_{\rho}(u_h) \subset X_h$. Thus with $w_h=u_{\epsilon,h}-v_{\epsilon,h}$, we obtain
\begin{align*}
\epsilon ||\Delta w_h||^2_0 - \int_{\Omega} (\det D^2 u_{\epsilon,h} - \det D^2 v_{\epsilon,h}) w_h \ud x = 0.
\end{align*}
Thus by the mean value theorem 
we have
\begin{align*}
\epsilon ||\Delta w_h||^2_0 - \int_{\Omega} [ \div ((\cof  (t D^2 u_{\epsilon,h} + (1-t) D^2 v_{\epsilon,h}) ) D w_h(x))] \cdot D w_h (x) \ud x & = 0 \\
\epsilon ||\Delta w_h||^2_0 + \int_{\Omega}  [\cof  (t D^2 u_{\epsilon,h} + (1-t) D^2 v_{\epsilon,h}) D w_h(x)] \cdot D w_h (x) \ud x & = 0.
\end{align*}
Using \eqref{p-d-K}, we obtain
\begin{align*}
0& =\epsilon ||\Delta w_h||^2_0 + \int_{\Omega}  [(\cof  D^2 w_h(x)) D w_h(x)] \cdot D w_h (x) \ud x  \\
& \geq \epsilon ||\Delta w_h||^2_0 + m |w_h|_1^2.
\end{align*}
Thus $|w_h|_1=0$ and since $w_h=0$ on $\partial \Omega$, we obtain $w_h=0$, the uniqueness of the discrete solution and the proof of the claim.


\section{Numerical results} \label{numerical}
The iterative method \eqref{newtonvm} depends on a parameter $\epsilon$ which has to be carefully chosen. As an alternative we present numerical results for a parameter independent iterative method. 
The latter is difficult to analyze and one can only expect a linear convergence rate. The numerical results are presented in order to illustrate some open problems in the numerical resolution of Monge-Amp\`ere equations. 

\subsection{Spline element method} \label{splineS}
The spline element method has been described in \cite{Awanou2003, Awanou2005a, Awanou2006, Baramidze2006, Hu2007}
under different names and more recently in \cite{Awanou2008}.
It can be described as a conforming discretization implementation with Lagrange multipliers. We first outline the main steps
of the method, discuss its advantages and possible disadvantages. We then give more details of this approach but refer to
the above references for explicit formulas. 

First, start with a representation of a piecewise discontinuous polynomial as a vector in $\mathbb{R}^N$, for some integer $N>0$.
Then express boundary conditions and constraints including global continuity or smoothness conditions as linear relations.
In our work, we use the Bernstein basis representation, \cite{Awanou2003,Awanou2008} which is very convenient to express smoothness conditions
and very popular in computer aided geometric design. Hence the term ``spline'' in the name of the method.
Splines are piecewise polynomials
with smoothness properties.
One then write a discrete version of the equation along with a discrete version of the spaces of trial and test functions.
The boundary conditions and constraints are enforced using Lagrange multipliers.
We are lead to saddle point problems which are solved by an augmented Lagrangian algorithm (sequences of linear equations with size $N \times N$).
The approach here should be contrasted with other approaches where Lagrange multipliers are introduced before discretization,
i.e. in \cite{Babuska73} or the discontinuous Galerkin methods.
 
The spline
element method, stands out as a robust,
flexible, efficient and accurate method.
It can be applied to a wide range of PDEs in science and engineering in both two and
three dimensions;
constraints and smoothness are enforced exactly and there is no need to implement
basis functions with the required properties; it is particularly suitable for fourth
order PDEs;
no inf-sup condition are needed to approximate Lagrange multipliers which arise
due to the constraints, e.g. the pressure term in the Navier-Stokes equations;
one gets in a single implementation approximations of variable order.
Other advantages of the method include the flexibility
of using polynomials of different degrees on different elements \cite{Hu2007}, the facility 
of implementing boundary conditions and the simplicity
of a posteriori error estimates since the method is conforming for many problems.
A possible disadvantage of this approach is the high number of degrees of freedom and the need to solve saddle point problems.

For illustration, we consider a general variational problem:  Find $u \in W$ such that
\begin{equation}
a(u,v) = \<l,v\> \quad \text{for all} \ v \in V, \label{var3}
\end{equation}
where $W$ and $V$ are respectively the space of trial and test functions.
We will assume that
the form $l$ is  bounded and linear and $a$ is a continuous 
mapping in some sense on $W \times V$ which is linear in the argument $v$. 

Let $W_h$ and $V_h$ be conforming subspaces of $W$ and $V$ respectively. We can write
$$
W_h =\{ c \in \mathbf{R}^N, R c =G \}, \ V_h = \{ c \in \mathbf{R}^N, R c =0 \},
$$
for a suitable vector $G$ and $R$ a suitable matrix which encodes the constraints
on the solution, e.g. smoothness and boundary conditions. 

The condition $a(u,v)=\<l,v\>$ for all $v \in V$ translates to
$$
K(c)d = L^T d \quad \forall d \in V_h, \ \text{that is for all} \ d \ \text{with} \ R d =0,
$$
for a suitable matrix $K(c)$ which depends on $c$ and $L$ is a vector of coefficients 
associated to the linear form $l$. If for 
example $\<l,v\>= \int_{\Omega} f v$, then $L^T d=d^T M F$ where $M$ is a mass matrix and $F$ a 
vector of coefficients associated to the spline interpolant of $f$. In the linear case
$K(c)$ can be written $c^T K$.

Introducing a Lagrange multiplier $\lambda$, 
the functional 
$$
K(c)d - L^Td + \lambda^T R d,
$$
vanishes identically on $V_h$. The stronger condition
$$
K(c) + \lambda^T R = L^T,
$$ 
along with the side condition $R c =G$ are the discrete equations to be solved.

By a slight abuse of notation, after linearization by Newton's method, 
the above nonlinear equation leads to solving systems of type
$$
c^T K  + \lambda^T R = L^T.
$$ 
The approximation $c$ of $u \in W$ thus is a limit of a sequence of
solutions of systems of type
\begin{align*}
\begin{split}
\left[\begin{array}{cc}K^T & R^T  \\ R & 0 \\
\end{array} \right]
\left[
\begin{array}{c} c \\ \lambda 
\end{array}\right] = \left[
\begin{array}{c} L\\G
\end{array}\right].
\end{split}
\end{align*}
It is therefore enough to consider the linear case.
If we assume for simplicity that $V=W$ and that the form $a$ is bilinear, symmetric, continuous and $V$-elliptic, existence of a discrete solution
follows from Lax-Milgram lemma. On the other hand, the ellipticity assures uniqueness of the component $c$ which can 
be retrieved by a least squares solution of the above system \cite{Awanou2003}. 
The Lagrange multiplier $\lambda$ may not be unique.
To avoid 
systems of large size, a variant of the augmented Lagrangian algorithm is used. 
For this, we consider the sequence of problems
\begin{align}
\begin{split} 
\left(\begin{array}{cc}K^T & R^T \\ R & -\mu M
\end{array} \right) \left[
\begin{array}{c} \mathbf{c}^{(l+1)} \\ \mathbf{\lambda}^{(l+1)}
\end{array}\right] = \left[
\begin{array}{c} L\\ G-\mu M\lambda^{(l)}
\end{array}\right],  \label{iter1}
\end{split}
\end{align}
where $\lambda^{(0)}$ is a suitable initial guess for example $\lambda^{(0)}=0$,
$M$ is a suitable matrix and $\mu>0$ is a small parameter taken in practice in the 
order of $10^{-5}$. It is possible to solve for $\mathbf{c}^{(l+1)}$ in terms of
$\mathbf{c}^{(l)}$.
A uniform convergence rate in $\mu$ for this algorithm was shown in \cite{Awanou2005}.


\subsection{Subharmonicity preserving iterations}
We also give numerical results for the following iterative method. Given an initial guess $u^0_h \in V_h$ with $u_h^{0}  = g_h$ $ \text{on} \, \partial \Omega$, find $u^{k+1}_h \in V_h$ such that $u_h^{k+1}  = g_h$ $\text{on} \, \partial \Omega$ and
\begin{equation} \label{subh-p}
\int_{\Omega} D u^{k+1}_h \cdot D v_h \ud x = - \int_{\Omega}  ((\Delta u^k_h)^n + n^n (f-\det  D^2 u^k_h ))^{\frac{1}{n}} v_h \ud x,  \forall v_h \in V_h \cap H_0^1(\Omega).
\end{equation}
The iterative method \eqref{subh-p} is the discrete analogue of the iterative method
\begin{equation}
\Delta u^{k+1} = ((\Delta u^k)^n + n^n (f-\det  D^2 u^k))^{\frac{1}{n}} \, \text{in} \, \Omega,  u^{k+1}=g \, \text{on} \,  \partial \Omega. \label{power}
\end{equation}
Since 
\begin{equation} \label{ag-ineq}
\det D^2 u^k \leq \frac{1}{n^n} (\Delta u^k)^n,
\end{equation}
it follows from \eqref{power} that $\Delta u^{k+1} \geq 0$. Hence, starting with an initial guess $u^0$ with $\Delta u^{0} \geq 0$, \eqref{power} preserves subharmonicity.
At the formal limit $\det D^2 u =f \geq 0$. Thus convexity is enforced for the two dimensional problem (at the continuous level). 
The iterative method \eqref{power} generalizes the method
\begin{equation}
\Delta u_{k+1} = ((\Delta u_k)^2 + 2(f-\det  D^2 u_k))^{\frac{1}{2}} \, \text{in} \, \Omega, u_{k+1} =g \, \text{on} \, \partial \Omega. \label{power2D}
\end{equation}
proposed in \cite{Benamou2010}.  In \cite{Oberman2010a,Oberman2010b}, the following generalization was proposed
\begin{equation}
\Delta u_{k+1} = ((\Delta u_k)^n + n!(f-\det  D^2 u_k))^{\frac{1}{2}} \, \text{in} \, \Omega, u_{k+1} =g \, \text{on} \, \partial \Omega.  \label{faux}
\end{equation}
It is clear that \eqref{power} and \eqref{faux} are different. Moreover, \eqref{power} is better since \eqref{faux} may not converge for a class of smooth functions as we now show. For $n>2$,  the method  \eqref{faux} can only converge for solutions of \eqref{m1} which also satisfies $(\Delta u)^2 = (\Delta u)^n$. Thus even for smooth solutions, the generalization we propose is better.

Let $a$ be such that $0 < a \leq n^n$. Then by \eqref{ag-ineq}, we have
$$
a \det D^2 v \leq (\Delta v)^n,
$$
and we can equally consider 
the iterative method
\begin{equation*}
\Delta u^{k+1} = ((\Delta u^k)^n + a (f-\det  D^2 u^k))^{\frac{1}{n}}, 
\end{equation*}
For $n=2$ and $a=2$ we get the one used in \cite{Benamou2010}. It will be referred to as the BFO iterative method.

In three dimension, we can also consider
\begin{equation}
\Delta u_{k+1} = ((\Delta u_k)^3 + 9(f-\det  D^2 u_k))^{\frac{1}{3}},  \label{powerm}
\end{equation}
corresponding to $a=9$. 

However, the formulation \eqref{power}, which shall henceforth be referred to as natural iterative method, appears to be better and this is supported numerically by a 2D example.

\subsection{Initial guess for the iterative methods}

The initial guess for the subharmonicity preserving  iterations is taken as the spline approximation of the solution of the Poisson equation
$\Delta u =  n^n f^{1/n}, \, n=2,3 \, \text{in} \, \Omega, \, u=g \, \text{on} \, \partial \Omega$. The initial guess for the discrete vanishing moment methodology is taken as the spline approximation of the biharmonic regularization
of a Poisson equation, 
$-\epsilon \Delta^2 u + \Delta u =  n f^{1/n} \, n=2,3 \,  \text{in}  \, \Omega, \, u=g, \Delta u =\epsilon^2 \, \text{on} \, \partial \Omega$.

\subsection{Two dimensional computational results}

The computational domain is the unit square $[0,1]^2$
which is first divided into squares of side length $h$. Then each square is 
divided into two triangles by the diagonal with negative slope. We recall that $d$ refers to the local degree of the piecewise polynomial
used.

For $g=0$, equation \eqref{m1} admits both a concave solution and a convex solution. Approximating concave solutions can be done by
either changing the initial guess or the structure of
the approximations.
\begin{enumerate}
\item Newton's method: initial guess  $\pm u_0$,
\item Iterative method \eqref{power2D}: $u_{k+1} = \pm \sqrt{(\Delta u_k)^2 + 2(f-\det  D^2 u_k)}$. 
\end{enumerate}

We consider the following test cases

Test 1: A smooth solution $u(x,y)=e^{(x^2+y^2)/2}$ so that 
$f(x,y)=(1+x^2+y^2)e^{(x^2+y^2)}$ and $g(x,y)=e^{(x^2+y^2)/2}$ on $\partial \Omega$.

The subharmonicity preserving iterations can be used with $C^0$ approximations. We compare the performance 
of \eqref{power} and \eqref{power2D} with Lagrange finite elements in Tables \ref{bfo-2d-1} and \ref{nat-2d-1}.


\begin{table}
\begin{tabular}{|c|c|c|c|c|} \hline
$h$  & $L^{2}$ norm & rate & $H^{1}$ norm  & rate \\ \hline
$1/2^1$   &  1.3558 $10^{-5}$& &  1.1212 $10^{-4}$ &     \\ \hline
$1/2^2$  &   9.2704 $10^{-7}$&3.87  & 5.5654 $10^{-6}$ & 4.33   \\ \hline
$1/2^3$  &   5.8359 $10^{-8}$&3.99 &  3.0329 $10^{-7}$ & 4.20   \\ \hline
$1/2^4$  &   3.6861 $10^{-9}$&3.98 &  1.8180 $10^{-8}$ &  4.06    \\ \hline
\end{tabular}
\bigskip
\caption{BFO method for Test 1, Lagrange elements $d=5$} \label{bfo-2d-1}
\end{table}
\begin{table}
\begin{tabular}{|c|c|c|c|c|} \hline
$h$  &   $L^{2}$ norm &  rate & $H^{1}$ norm & rate \\ \hline
$1/2^1$       &   3.4383 $10^{-6}$& &  8.8363  $10^{-5}$ &  \\ \hline
$1/2^2$     &   1.1022 $10^{-7}$&4.96 &  3.1305 $10^{-6}$  & 4.82  \\ \hline
$1/2^3$    &   7.5096 $10^{-9}$ &3.87 &  1.0762 $10^{-7}$   &4.86   \\ \hline
$1/2^4$    &   4.9561 $10^{-10}$& 3.92 & 4.1682 $10^{-9}$  &4.69   \\ \hline
\end{tabular}
\bigskip
\caption{Natural iterative method for Test 1, Lagrange elements $d=5$} \label{nat-2d-1}
\end{table}

Test 2: A solution not in $H^2(\Omega)$, $u(x,y)=-\sqrt{2-x^2-y^2}$ 
so that $f(x,y)=2/(2-x^2-y^2)^2$ and $g(x,y)=-\sqrt{2-x^2-y^2}$ on $\partial \Omega$.


For this non smooth solution, it is essential to adapt the choice of the parameter $\epsilon$ in the discrete vanishing methodology as a function of $h$, c.f. Table \ref{vm-var}. This point was already made in \cite{Feng2009}. We also recall that for this problem, Newton's method diverges and note that the subharmonicity preserving iterations are robust, c.f. Table \ref{sub-var}.

\begin{table}
\begin{minipage}{2in} 
\begin{tabular}{|c|c|c|} \hline
$h$  & $L^{2}$ norm & $H^{1}$ norm  \\ \hline
$1/2^1$   & 7.6680$10^{-3}$  & 7.4491$10^{-2}$       \\ \hline
$1/2^2$  &  1.4536$10^{-3}$ &  3.9244$10^{-2}$    \\ \hline
$1/2^3$  &  9.8727$10^{-3}$ &  2.5112$10^{-1}$    \\ \hline
$1/2^4$  &  5.6819$10^{-3}$ &  2.4927$10^{-1}$       \\ \hline
$1/2^5$  &  1.9830 $10^{+4}$ &  1.1812 $10^{+6}$      \\ \hline
\end{tabular}
\end{minipage}
\hspace{1cm}
\begin{minipage}{2in} 
\begin{tabular}{|c|c|c|} \hline
$h$  & $L^{2}$ norm & $H^{1}$ norm  \\ \hline
$1/2^1$   & 7.8254$10^{-3}$ &  9.3184$10^{-2}$        \\ \hline
$1/2^2$  &  1.0646$10^{-2}$ & 9.5201$10^{-2}$      \\ \hline
$1/2^3$  &  1.1306$10^{-2}$ & 9.6154$10^{-2}$     \\ \hline
$1/2^4$  & 1.1500$10^{-2}$ & 9.1336$10^{-2}$         \\ \hline
$1/2^5$  &  1.1625$10^{-2}$ & 8.7785$10^{-2}$      \\ \hline
 $1/2^6$ &  1.1681$10^{-2}$ & 8.5632$10^{-2}$                     \\ \hline
\end{tabular}
\end{minipage}
\bigskip
\caption{ Vanishing moment Test 2 $\epsilon=10^{-3}$ and $\epsilon=10^{-2}$, $d=5$} \label{vm-var}
\end{table}

\begin{table}
\begin{minipage}{2in} %
\begin{tabular}{|c|c|c|} \hline
$h$  & $L^{2}$ norm & $H^{1}$ norm  \\ \hline
$1/2^1$   &  2.1954$10^{-2}$ &  1.6409$10^{-1}$      \\ \hline
$1/2^2$  &   3.6097$10^{-3}$  & 6.1405$10^{-2}$     \\ \hline
$1/2^3$  &   1.0685$10^{-3}$ &  4.0978$10^{-2}$     \\ \hline
$1/2^4$  &   5.0838$10^{-3}$ &  2.8048$10^{-1}$       \\ \hline
$1/2^5$  &   2.5797$10^{+3}$ &  2.2688$10^{+5}$      \\ \hline
$1/2^6$  &   1.8452$10^{+4}$ &  3.5922$10^{+6}$      \\ \hline
\end{tabular}
\end{minipage}
\hspace{1cm}
\begin{minipage}{2in} 
\begin{tabular}{|c|c|c|c|} \hline
$h$  &  $n_{\text{it}}$ & $L^{2}$ norm & $H^{1}$ norm  \\ \hline
$1/2^1$   & 50    &   2.3921 $10^{-1}$ &  1.1900     \\ \hline
$1/2^2$  &  159   &   1.2585 $10^{-1}$ &  7.1292 $10^{-1}$     \\ \hline
$1/2^3$  &  151   &   1.0341 $10^{-1}$ &  6.4299 $10^{-1}$      \\ \hline
$1/2^4$  & 160   &   9.6031 $10^{-2}$ &  6.2088 $10^{-1}$     \\ \hline
$1/2^5$  &     199   &   9.4551 $10^{-2}$ &  6.2453 $10^{-1}$   \\ \hline
$1/2^6$  & 8     &   1.6977 $10^{-2}$ &  2.2925 $10^{-1}$   \\ \hline
\end{tabular}
\end{minipage}
\bigskip
\caption{Newton's method and BFO iterative method for Test 2, $d=3$} \label{sub-var}
\end{table}

\subsection{Three dimensional computational results}

B$\ddot{\text{o}}$hmer \cite{Bohmer2008} and Feng and Neilan \cite{Feng2009} have discussed the possibility of using $C^1$ finite elements in  three dimensions but no numerical evidence was given. This can be addressed with the spline element method. We used two computational domains both on the unit cube $[0,1]^3$
which is first divided into six tetrahedra (Domain 1) or twelve tetrahedra (Domain 2) forming a tetrahedral partition $\mathcal{T}_1$. 
This partition is uniformly refined
following a strategy introduced in \cite{Awanou2003} similar to the one of \cite{Ong94}
resulting in successive level of refinements $\mathcal{T}_k$, $k=2,3,\ldots$

We consider the following test cases

Test 3: $u(x,y,z)=e^{(x^2+y^2+z^2)/3}$ so that
$f(x,y,z)=8/81(3+2(x^2+y^2+z^2))e^{(x^2+y^2+z^2)}$ and $g(x,y,z)=e^{(x^2+y^2+z^2)/3}$ on 
$\partial \Omega$. 

Since the solution is smooth, it is enough to use Newton's method, c.f. Tables \ref{newton-3d-1} and \ref{newton-3d-2}. To emphasize this point, we numerically show in Table \ref{num-cvg-vm} the convergence as $\epsilon \to 0$ of the solution of \eqref{vm} to the solution of \eqref{m1h}.



\begin{table}
\begin{tabular}{|c|c|c|c|} \hline
d & $L^{2}$ norm & $H^{1}$ norm &  $H^{2}$ norm \\ \hline
3 & 1.2338 $10^{-2}$ &  7.6984 $10^{-2}$ &    4.4411 $10^{-1}$ \\ \hline
4 & 1.6289 $10^{-3}$ &  1.4719 $10^{-2}$ &    1.3983 $10^{-1}$ \\ \hline
5 & 1.5333 $10^{-3}$ &  8.7312 $10^{-3}$ &    6.0412 $10^{-2}$\\ \hline
6 & 1.2324 $10^{-4}$ &  9.7171 $10^{-4}$ &    1.0584 $10^{-2}$ \\ \hline
Rate & 0.18 $0.25^{d-1}$    & 4.58 $0.25^{d}$   & 59.96 $0.3^{d+1}$   \\ \hline
\end{tabular}
\bigskip
\caption{Newton's method Test 3, Domain 1 on $ \mathcal{I}_1 $} \label{newton-3d-1}
\end{table}

\begin{table}
\begin{tabular}{|c|c|c|c|c|} \hline
d & $L^{2}$ norm & $H^{1}$ norm &  $H^{2}$ norm \\ \hline
3 & 3.1739 $10^{-3}$ &  2.3005 $10^{-2}$ &   2.4496 $10^{-1}$ \\ \hline
4 & 3.2786 $10^{-4}$ &  3.5626 $10^{-3}$ &   5.2079 $10^{-2}$  \\ \hline
5 & 2.4027 $10^{-5}$ &  3.9210 $10^{-4}$ &   8.8868 $10^{-3}$  \\ \hline
6 & 1.3821 $10^{-6}$ &  2.2369 $10^{-5}$ &   6.0918 $10^{-4}$  \\ \hline
Rate & 0.65 $0.075^{d-1}$    &  28.96 $0.1^d$  & 849.85 $0.14^{d+1}$   \\ \hline
\end{tabular}
\bigskip
\caption{Newton's method Test 3, Domain 1 on $\mathcal{T}_2$}  \label{newton-3d-2}
\end{table}

\begin{table}\label{3Drobust}
\begin{tabular}{|c|c|c|c|} \hline
$\epsilon$  & $L^{2}$ norm & $H^{1}$ norm &  $H^{2}$ norm \\ \hline
$10^{-1}$ & 6.6870 $10^{-2}$ &  3.9292 $10^{-1}$ &  2.8852 \\ \hline
$10^{-2}$  & 1.8832 $10^{-2}$ &  1.3137 $10^{-1}$ &  1.5882  \\ \hline
$10^{-3}$  & 2.4237 $10^{-3}$ &  2.5273 $10^{-2}$ &  5.3206 $10^{-1}$   \\ \hline
$10^{-4}$  & 2.5661 $10^{-4}$ &  3.2633 $10^{-3}$ &  7.9936 $10^{-2}$  \\ \hline
$10^{-5}$  & 3.1058 $10^{-5}$ &  5.0367 $10^{-4}$ &  1.2543 $10^{-2}$ \\ \hline
$10^{-6}$  &  2.3519 $10^{-5}$ &  3.9165 $10^{-4}$ &  8.9744 $10^{-3}$ \\ \hline
$10^{-7}$  & 2.3964 $10^{-5}$ &  3.9193 $10^{-4}$ &  8.8921 $10^{-3}$ \\ \hline
$10^{-10}$  & 2.4027 $10^{-5}$ &  3.9210 $10^{-4}$ &  8.8868 $10^{-3}$ \\ \hline  
0  & 2.4027 $10^{-5}$ &  3.9210 $10^{-4}$ &   8.8868 $10^{-3}$   \\ \hline
\end{tabular}    
\bigskip
\caption{3D numerical robustness Test 5, Domain 1 on $\mathcal{T}_2$, $d=5$} \label{num-cvg-vm}
\end{table}

\begin{figure}[tbp]
\begin{center}
\includegraphics[angle=0, height=7cm]{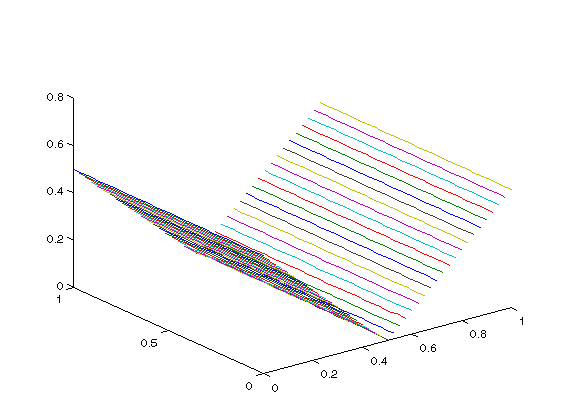}
\end{center}
\bigskip
\caption{Method \eqref{powerm} for Test 4 on Domain 2 and $ \mathcal{I}_3$, $d=5, r=1$, plane $z=0$.} \label{fig-0}
\end{figure}

Test 4: $f(x,y,z)=0$ and $g(x,y,z)=|x-1/2|$.  For the degenerate case of this test, we did not capture the convexity of the discrete solution by discretizing \eqref{powerm} with the standard finite difference method. Surprisingly, with $C^1$ splines, we were able to capture a $C^1$ function which appears to approximate well the solution, c.f. Figure \ref{fig-0}.

For non smooth solutions in three dimension, it would be better to have iterative methods which preserve explicitly convexity. 
In some cases, i.e. for $f=1$ and $g=0$, we were able to capture the correct solution. Our methods did not work for the 3D analogue of Test 2. 

\providecommand{\bysame}{\leavevmode\hbox to3em{\hrulefill}\thinspace}
\providecommand{\MR}{\relax\ifhmode\unskip\space\fi MR }
\providecommand{\MRhref}[2]{%
  \href{http://www.ams.org/mathscinet-getitem?mr=#1}{#2}
}
\providecommand{\href}[2]{#2}
\begin{thebibliography}{10}

\bibitem{Aleksandrov1961}
A.~D. Aleksandrov, \emph{Certain estimates for the {D}irichlet problem}, Soviet
  Math. Dokl. \textbf{1} (1961), 1151--1154.

\bibitem{Aubin1982}
Thierry Aubin, \emph{Nonlinear analysis on manifolds. {M}onge-{A}mp\`ere
  equations}, Grundlehren der Mathematischen Wissenschaften [Fundamental
  Principles of Mathematical Sciences], vol. 252, Springer-Verlag, New York,
  1982.

\bibitem{Avellaneda1994}
M.~Avellaneda and A.~Paras, \emph{Dynamic hedging portfolios for derivative
  securities in the presence of large transaction costs}, Appl. Math. Finance
  \textbf{1} (1994), 165--194.

\bibitem{Awanou2003}
G.~Awanou, \emph{Energy methods in 3{D} spline approximations of the
  {N}avier-{S}tokes equations}, Ph.D. Dissertation, University of Georgia,
  Athens, Ga, 2003.

\bibitem{Awanou2005}
G.~M. Awanou and M.~J. Lai, \emph{On convergence rate of the augmented
  {L}agrangian algorithm for nonsymmetric saddle point problems}, Appl. Numer.
  Math. \textbf{54} (2005), no.~2, 122--134.

\bibitem{Awanou2008}
Gerard Awanou, \emph{Robustness of a spline element method with constraints},
  J. Sci. Comput. \textbf{36} (2008), no.~3, 421--432.

\bibitem{Awanou2005a}
Gerard Awanou and Ming-Jun Lai, \emph{Trivariate spline approximations of 3{D}
  {N}avier-{S}tokes equations}, Math. Comp. \textbf{74} (2005), no.~250,
  585--601 (electronic).

\bibitem{Awanou2006}
Gerard Awanou, Ming-Jun Lai, and Paul Wenston, \emph{The multivariate spline
  method for scattered data fitting and numerical solution of partial
  differential equations}, Wavelets and splines: {A}thens 2005, Mod. Methods
  Math., Nashboro Press, Brentwood, TN, 2006, pp.~24--74.

\bibitem{Babuska73}
Ivo Babu{\v{s}}ka, \emph{The finite element method with {L}agrangian
  multipliers}, Numer. Math. \textbf{20} (1972/73), 179--192.

\bibitem{Bakelman65}
I.~Ja. Bakelman, \emph{Geometricheskie metody resheniya ellipticheskikh
  uravnenii}, Izdat. ``Nauka'', Moscow, 1965.

\bibitem{Bakelman1983}
Ilya~J. Bakelman, \emph{Variational problems and elliptic {M}onge-{A}mp\`ere
  equations}, J. Differential Geom. \textbf{18} (1983), no.~4, 669--699 (1984).

\bibitem{Bakelman1994}
\bysame, \emph{Convex analysis and nonlinear geometric elliptic equations},
  Springer-Verlag, Berlin, 1994, With an obituary for the author by William
  Rundell, Edited by Steven D. Taliaferro.

\bibitem{Baramidze2006}
Victoria Baramidze and Ming-Jun Lai, \emph{Spherical spline solution to a {PDE}
  on the sphere}, Wavelets and splines: {A}thens 2005, Mod. Methods Math.,
  Nashboro Press, Brentwood, TN, 2006, pp.~75--92.

\bibitem{Barles1991}
G.~Barles and P.~E. Souganidis, \emph{Convergence of approximation schemes for
  fully nonlinear second order equations}, Asymptotic Anal. \textbf{4} (1991),
  no.~3, 271--283.

\bibitem{Benamou2000}
Jean-David Benamou and Yann Brenier, \emph{A computational fluid mechanics
  solution to the {M}onge-{K}antorovich mass transfer problem}, Numer. Math.
  \textbf{84} (2000), no.~3, 375--393.

\bibitem{Benamou2010}
Jean-David Benamou, Brittany Froese, and Adam Oberman, \emph{Two numerical
  methods for the elliptic {M}onge-{A}mp\`ere equation}, 2010.

\bibitem{Belgacem2006}
M.~Bouchiba and F.~Ben Belgacem, \emph{Numerical solution of {M}onge-{A}mpere
  equation}, Math. Balkanica (N.S.) \textbf{20} (2006), no.~3-4, 369--378.

\bibitem{Brenner2010a}
S.~Brenner and M.~Neilan, \emph{Finite element approximations of the three
  dimensional {M}onge-{A}mpere equation}, Submitted, 2010.

\bibitem{Caffarelli1995}
Luis~A. Caffarelli and Xavier Cabr{\'e}, \emph{Fully nonlinear elliptic
  equations}, American Mathematical Society Colloquium Publications, vol.~43,
  American Mathematical Society, Providence, RI, 1995.

\bibitem{Caffarelli1999}
Luis~A. Caffarelli, Sergey~A. Kochengin, and Vladimir~I. Oliker, \emph{On the
  numerical solution of the problem of reflector design with given far-field
  scattering data}, Monge {A}mp\`ere equation: applications to geometry and
  optimization ({D}eerfield {B}each, {FL}, 1997), Contemp. Math., vol. 226,
  Amer. Math. Soc., Providence, RI, 1999, pp.~13--32.

\bibitem{Carlen2004}
E.~A. Carlen and W.~Gangbo, \emph{Solution of a model {B}oltzmann equation via
  steepest descent in the 2-{W}asserstein metric}, Arch. Ration. Mech. Anal.
  \textbf{172} (2004), no.~1, 21--64.

\bibitem{Chen98}
Ya-Zhe Chen and Lan-Cheng Wu, \emph{Second order elliptic equations and
  elliptic systems}, Translations of Mathematical Monographs, vol. 174,
  American Mathematical Society, Providence, RI, 1998, Translated from the 1991
  Chinese original by Bei Hu.

\bibitem{ChengYau77}
Shiu~Yuen Cheng and Shing~Tung Yau, \emph{On the regularity of the
  {M}onge-{A}mp\`ere equation {${\rm det}(\partial ^{2}u/\partial x_{i}\partial
  sx_{j})=F(x,u)$}}, Comm. Pure Appl. Math. \textbf{30} (1977), no.~1, 41--68.
  \MR{MR0437805 (55 \#10727)}

\bibitem{Courant1989}
R.~Courant and D.~Hilbert, \emph{Methods of mathematical physics. {V}ol. {II}},
  Wiley Classics Library, John Wiley \& Sons Inc., New York, 1989, Partial
  differential equations, Reprint of the 1962 original, A Wiley-Interscience
  Publication.

\bibitem{Crandall1992}
Michael~G. Crandall, Hitoshi Ishii, and Pierre-Louis Lions, \emph{User's guide
  to viscosity solutions of second order partial differential equations}, Bull.
  Amer. Math. Soc. (N.S.) \textbf{27} (1992), no.~1, 1--67.

\bibitem{Cullen2001}
Mike Cullen and Wilfrid Gangbo, \emph{A variational approach for the
  2-dimensional semi-geostrophic shallow water equations}, Arch. Ration. Mech.
  Anal. \textbf{156} (2001), no.~3, 241--273.

\bibitem{Dean2006}
E.~J. Dean and R.~Glowinski, \emph{Numerical methods for fully nonlinear
  elliptic equations of the {M}onge-{A}mp\`ere type}, Comput. Methods Appl.
  Mech. Engrg. \textbf{195} (2006), no.~13-16, 1344--1386.

\bibitem{Dean2003}
Edward~J. Dean and Roland Glowinski, \emph{Numerical solution of the
  two-dimensional elliptic {M}onge-{A}mp\`ere equation with {D}irichlet
  boundary conditions: an augmented {L}agrangian approach}, C. R. Math. Acad.
  Sci. Paris \textbf{336} (2003), no.~9, 779--784.

\bibitem{Dean2004}
\bysame, \emph{Numerical solution of the two-dimensional elliptic
  {M}onge-{A}mp\`ere equation with {D}irichlet boundary conditions: a
  least-squares approach}, C. R. Math. Acad. Sci. Paris \textbf{339} (2004),
  no.~12, 887--892.

\bibitem{Dyer2003}
Bradley~W. Dyer and Don Hong, \emph{Algorithm for optimal triangulations in
  scattered data representation and implementation}, J. Comput. Anal. Appl.
  \textbf{5} (2003), no.~1, 25--43, Approximation theory and wavelets (Austin,
  TX, 1999).

\bibitem{Evans1998}
Lawrence~C. Evans, \emph{Partial differential equations}, Graduate Studies in
  Mathematics, vol.~19, American Mathematical Society, Providence, RI, 1998.

\bibitem{Feng2009}
X.~Feng and M.~Neilan, \emph{Analysis of {G}alerkin methods for the fully
  nonlinear {M}onge-{A}mp\`ere equation}, Preprint, 2009.

\bibitem{Feng2009a}
\bysame, \emph{Error analysis for mixed finite element approximations of the
  fully nonlinear monge-amp\`ere equation based on the vanishing moment
  method}, SIAM J. Numer. Anal. \textbf{47} (2009), no.~2, 1226--1250.

\bibitem{Feng2009b}
\bysame, \emph{Vanishing moment method and moment solutions for second order
  fully nonlinear partial differential equations}, J. Sci. Comput. \textbf{38}
  (2009), no.~1, 74--98.

\bibitem{Oberman2010a}
B.D. Froese and A.M. Oberman, \emph{Convergent finite difference solvers for
  viscosity solutions of the elliptic {M}onge-{A}mpere equation in dimensions
  two and higher}, Submitted, 2010.

\bibitem{Oberman2010b}
\bysame, \emph{Fast convergent finite difference methods for the elliptic
  {M}onge-{A}mpere equation in dimensions two and higher}, Submitted, 2010.

\bibitem{Gangbo2004}
W.~Gangbo., \emph{An introduction to the mass transportation theory and its
  applications}, Carnegie Mellon 2004 Summer Institute, 2004.

\bibitem{Gangbo2000}
Wilfrid Gangbo and Roberto Van~der Putten, \emph{Uniqueness of equilibrium
  configurations in solid crystals}, SIAM J. Math. Anal. \textbf{32} (2000),
  no.~3, 465--492 (electronic).

\bibitem{Guti'errez2001}
Cristian~E. Guti{\'e}rrez, \emph{The {M}onge-{A}mp\`ere equation}, Progress in
  Nonlinear Differential Equations and their Applications, 44, Birkh\"auser
  Boston Inc., Boston, MA, 2001.

\bibitem{Hartman66}
Philip Hartman and Guido Stampacchia, \emph{On some non-linear elliptic
  differential-functional equations}, Acta Math. \textbf{115} (1966), 271--310.

\bibitem{Hu2007}
Xian-Liang Hu, Dan-Fu Han, and Ming-Jun Lai, \emph{Bivariate splines of various
  degrees for numerical solution of partial differential equations}, SIAM J.
  Sci. Comput. \textbf{29} (2007), no.~3, 1338--1354 (electronic).

\bibitem{Kochengin1998}
Sergey~A. Kochengin and Vladimir~I. Oliker, \emph{Determination of reflector
  surfaces from near-field scattering data. {II}. {N}umerical solution}, Numer.
  Math. \textbf{79} (1998), no.~4, 553--568.

\bibitem{Lai2007}
Ming-Jun Lai and Larry~L. Schumaker, \emph{Spline functions on triangulations},
  Encyclopedia of Mathematics and its Applications, vol. 110, Cambridge
  University Press, Cambridge, 2007.

\bibitem{Lindenstrauss94}
Joram Lindenstrauss, Lawrence~C. Evans, Adrien Douady, Aner Shalev, and
  Nicholas Pippenger, \emph{Fields {M}edals and {N}evanlinna {P}rize presented
  at {ICM}-94 in {Z}\"urich}, Notices Amer. Math. Soc. \textbf{41} (1994),
  no.~9, 1103--1111.

\bibitem{Loeper2005}
Gr{\'e}goire Loeper and Francesca Rapetti, \emph{Numerical solution of the
  {M}onge-{A}mp\`ere equation by a {N}ewton's algorithm}, C. R. Math. Acad.
  Sci. Paris \textbf{340} (2005), no.~4, 319--324.

\bibitem{Mohammadi2007}
Bijan Mohammadi, \emph{Optimal transport, shape optimization and global
  minimization}, C. R. Math. Acad. Sci. Paris \textbf{344} (2007), no.~9,
  591--596.

\bibitem{Winther2001}
Trygve~K. Nilssen, Xue-Cheng Tai, and Ragnar Winther, \emph{A robust
  nonconforming {$H^2$}-element}, Math. Comp. \textbf{70} (2001), no.~234,
  489--505.

\bibitem{Oberman2008}
Adam~M. Oberman, \emph{Wide stencil finite difference schemes for the elliptic
  {M}onge-{A}mp\`ere equation and functions of the eigenvalues of the
  {H}essian}, Discrete Contin. Dyn. Syst. Ser. B \textbf{10} (2008), no.~1,
  221--238.

\bibitem{Ockendon2003}
John Ockendon, Sam Howison, Andrew Lacey, and Alexander Movchan, \emph{Applied
  partial differential equations}, revised ed., Oxford University Press,
  Oxford, 2003.

\bibitem{Oliker1988}
V.~I. Oliker and L.~D. Prussner, \emph{On the numerical solution of the
  equation {$(\partial\sp 2z/\partial x\sp 2)(\partial\sp 2z/\partial y\sp
  2)-((\partial\sp 2z/\partial x\partial y))\sp 2=f$} and its discretizations.
  {I}}, Numer. Math. \textbf{54} (1988), no.~3, 271--293.

\bibitem{Ong94}
Maria Elizabeth~G. Ong, \emph{Uniform refinement of a tetrahedron}, SIAM J.
  Sci. Comput. \textbf{15} (1994), no.~5, 1134--1144.

\bibitem{Brenner2010b}
M.~Neilan S.~C.~Brenner, T.~Gudi and L.Y. Sung, \emph{C0 penalty methods for
  the fully nonlinear {M}onge-{A}mpere equation}, Submitted, 2010.

\bibitem{Schumaker2009}
Larry~L. Schumaker, Tatyana Sorokina, and Andrew~J. Worsey, \emph{A {$C^1$}
  quadratic trivariate macro-element space defined over arbitrary tetrahedral
  partitions}, J. Approx. Theory \textbf{158} (2009), no.~1, 126--142.
  \MR{2523724 (2010i:65023)}

\bibitem{Swarztrauber84}
Paul~N. Swarztrauber, \emph{Fast {P}oisson solvers},  \textbf{24} (1984),
  319--370.

\bibitem{Tornberg03}
Anna-Karin Tornberg and Bj{\"o}rn Engquist, \emph{Regularization techniques for
  numerical approximation of {PDE}s with singularities}, J. Sci. Comput.
  \textbf{19} (2003), no.~1-3, 527--552, Special issue in honor of the sixtieth
  birthday of Stanley Osher.

\bibitem{Trudinger08}
Neil~S. Trudinger and Xu-Jia Wang, \emph{Boundary regularity for the
  {M}onge-{A}mp\`ere and affine maximal surface equations}, Ann. of Math. (2)
  \textbf{167} (2008), no.~3, 993--1028.

\bibitem{Tso1990}
Kaising Tso, \emph{On a real {M}onge-{A}mp\`ere functional}, Invent. Math.
  \textbf{101} (1990), no.~2, 425--448. \MR{1062970 (91i:35082)}

\bibitem{Westcott1983}
B.S. Westcott, \emph{Shaped reflector antenna design}, Letchworth, U.K., 1983.

\bibitem{Zheligovsky10}
V.~Zheligovsky, O.~Podvigina, and U.~Frisch, \emph{The {M}onge-{A}mp\`ere
  equation: various forms and numerical solution}, J. Comput. Phys.
  \textbf{229} (2010), no.~13, 5043--5061.

\end{thebibliography}


\begin{thebibliography}{10}
\providecommand{\url}[1]{{#1}}
\providecommand{\urlprefix}{URL }
\expandafter\ifx\csname urlstyle\endcsname\relax
  \providecommand{\doi}[1]{DOI~\discretionary{}{}{}#1}\else
  \providecommand{\doi}{DOI~\discretionary{}{}{}\begingroup
  \urlstyle{rm}\Url}\fi

\bibitem{Awanou2003}
Awanou, G.: Energy methods in 3{D} spline approximations of the
  {N}avier-{S}tokes equations.
\newblock Ph.D. Dissertation, University of Georgia. Athens, Ga (2003)

\bibitem{Awanou2008}
Awanou, G.: Robustness of a spline element method with constraints.
\newblock J. Sci. Comput. \textbf{36}(3), 421--432 (2008)

\bibitem{AwanouPseudo10}
Awanou, G.: Pseudo transient continuation and time marching methods for
  {M}onge-{A}mp\`ere type equations (2013).
\newblock {h}ttp://arxiv.org/pdf/1301.5891.pdf

\bibitem{Awanou-Std04}
Awanou, G.: Standard finite elements for the numerical resolution of the
  elliptic {M}onge-{A}mp\`ere equation: {A}leksandrov solutions (2013).
\newblock {h}ttp://arxiv.org/pdf/1310.4568v1.pdf

\bibitem{Awanou2005a}
Awanou, G., Lai, M.J.: Trivariate spline approximations of 3{D}
  {N}avier-{S}tokes equations.
\newblock Math. Comp. \textbf{74}(250), 585--601 (electronic) (2005)

\bibitem{Awanou2006}
Awanou, G., Lai, M.J., Wenston, P.: The multivariate spline method for
  scattered data fitting and numerical solution of partial differential
  equations.
\newblock In: Wavelets and splines: {A}thens 2005, Mod. Methods Math., pp.
  24--74. Nashboro Press, Brentwood, TN (2006)

\bibitem{Awanou2005}
Awanou, G.M., Lai, M.J.: On convergence rate of the augmented {L}agrangian
  algorithm for nonsymmetric saddle point problems.
\newblock Appl. Numer. Math. \textbf{54}(2), 122--134 (2005)

\bibitem{Babuska73}
Babu{\v{s}}ka, I.: The finite element method with {L}agrangian multipliers.
\newblock Numer. Math. \textbf{20}, 179--192 (1972/73)

\bibitem{Baramidze2006}
Baramidze, V., Lai, M.J.: Spherical spline solution to a {PDE} on the sphere.
\newblock In: Wavelets and splines: {A}thens 2005, Mod. Methods Math., pp.
  75--92. Nashboro Press, Brentwood, TN (2006)

\bibitem{Benamou2010}
Benamou, J.D., Froese, B.D., Oberman, A.M.: Two numerical methods for the
  elliptic {M}onge-{A}mp\`ere equation.
\newblock M2AN Math. Model. Numer. Anal. \textbf{44}(4), 737--758 (2010)

\bibitem{Bohmer2008}
B{\"o}hmer, K.: On finite element methods for fully nonlinear elliptic
  equations of second order.
\newblock SIAM J. Numer. Anal. \textbf{46}(3), 1212--1249 (2008)

\bibitem{Bohmer2010}
Bohmer, K.: Numerical Methods for Nonlinear Elliptic Differential Equations: A
  Synopsis.
\newblock Oxford University Press, USA (2010)

\bibitem{Belgacem2006}
Bouchiba, M., Belgacem, F.B.: Numerical solution of {M}onge-{A}mpere equation.
\newblock Math. Balkanica (N.S.) \textbf{20}(3-4), 369--378 (2006)

\bibitem{Brenner2010b}
Brenner, S.C., Gudi, T., Neilan, M., Sung, L.Y.: {$C^0$} penalty methods for
  the fully nonlinear {M}onge-{A}mp\`ere equation.
\newblock Math. Comp. \textbf{80}(276), 1979--1995 (2011)

\bibitem{Brenner2010a}
Brenner, S.C., Neilan, M.: Finite element approximations of the three
  dimensional {M}onge-{A}mp\`ere equation.
\newblock ESAIM Math. Model. Numer. Anal. \textbf{46}(5), 979--1001 (2012)

\bibitem{Brenner02}
Brenner, S.C., Scott, L.R.: The mathematical theory of finite element methods,
  \emph{Texts in Applied Mathematics}, vol.~15, second edn.
\newblock Springer-Verlag, New York (2002)

\bibitem{Caffarelli1999}
Caffarelli, L.A., Kochengin, S.A., Oliker, V.I.: On the numerical solution of
  the problem of reflector design with given far-field scattering data.
\newblock In: Monge {A}mp\`ere equation: applications to geometry and
  optimization ({D}eerfield {B}each, {FL}, 1997), \emph{Contemp. Math.}, vol.
  226, pp. 13--32. Amer. Math. Soc., Providence, RI (1999)

\bibitem{Dean2003}
Dean, E.J., Glowinski, R.: Numerical solution of the two-dimensional elliptic
  {M}onge-{A}mp\`ere equation with {D}irichlet boundary conditions: an
  augmented {L}agrangian approach.
\newblock C. R. Math. Acad. Sci. Paris \textbf{336}(9), 779--784 (2003)

\bibitem{Dean2004}
Dean, E.J., Glowinski, R.: Numerical solution of the two-dimensional elliptic
  {M}onge-{A}mp\`ere equation with {D}irichlet boundary conditions: a
  least-squares approach.
\newblock C. R. Math. Acad. Sci. Paris \textbf{339}(12), 887--892 (2004)

\bibitem{Dean2006}
Dean, E.J., Glowinski, R.: Numerical methods for fully nonlinear elliptic
  equations of the {M}onge-{A}mp\`ere type.
\newblock Comput. Methods Appl. Mech. Engrg. \textbf{195}(13-16), 1344--1386
  (2006)

\bibitem{Feng2012}
Feng, X., Neilan, M.: Convergence of a fourth order singular perturbation of
  the $n$-dimensional radially symmetric {M}onge-{A}mpere equation.
\newblock Submitted

\bibitem{Feng2009a}
Feng, X., Neilan, M.: Error analysis for mixed finite element approximations of
  the fully nonlinear {M}onge-{A}mp\`ere equation based on the vanishing moment
  method.
\newblock SIAM J. Numer. Anal. \textbf{47}(2), 1226--1250 (2009)

\bibitem{Feng2009b}
Feng, X., Neilan, M.: Vanishing moment method and moment solutions for second
  order fully nonlinear partial differential equations.
\newblock J. Sci. Comput. \textbf{38}(1), 74--98 (2009)

\bibitem{Feng2009}
Feng, X., Neilan, M.: Analysis of {G}alerkin methods for the fully nonlinear
  {M}onge-{A}mp\`ere equation.
\newblock J. Sci. Comput. \textbf{47}(3), 303--327 (2011)

\bibitem{Oberman2010a}
Froese, B., Oberman, A.: Convergent finite difference solvers for viscosity
  solutions of the elliptic {M}onge-{A}mp\`ere equation in dimensions two and
  higher.
\newblock SIAM J. Numer. Anal. \textbf{49}(4), 1692--1714 (2011)

\bibitem{Oberman2010b}
Froese, B.D., Oberman, A.M.: Fast finite difference solvers for singular
  solutions of the elliptic {M}onge-{A}mp\`ere equation.
\newblock J. Comput. Phys. \textbf{230}(3), 818--834 (2011)

\bibitem{GlowinskiICIAM07}
Glowinski, R.: Numerical methods for fully nonlinear elliptic equations.
\newblock In: I{CIAM} 07---6th {I}nternational {C}ongress on {I}ndustrial and
  {A}pplied {M}athematics, pp. 155--192. Eur. Math. Soc., Z\"urich (2009)

\bibitem{Hu2007}
Hu, X.L., Han, D.F., Lai, M.J.: Bivariate splines of various degrees for
  numerical solution of partial differential equations.
\newblock SIAM J. Sci. Comput. \textbf{29}(3), 1338--1354 (electronic) (2007)

\bibitem{Kochengin1998}
Kochengin, S.A., Oliker, V.I.: Determination of reflector surfaces from
  near-field scattering data. {II}. {N}umerical solution.
\newblock Numer. Math. \textbf{79}(4), 553--568 (1998)

\bibitem{Lai98}
Lai, M.J., Schumaker, L.L.: On the approximation power of bivariate splines.
\newblock Adv. Comput. Math. \textbf{9}(3-4), 251--279 (1998)

\bibitem{Lai2007}
Lai, M.J., Schumaker, L.L.: Spline functions on triangulations,
  \emph{Encyclopedia of Mathematics and its Applications}, vol. 110.
\newblock Cambridge University Press, Cambridge (2007)

\bibitem{Lakkis11b}
Lakkis, O., Pryer, T.: A {F}inite {E}lement {M}ethod for {N}onlinear {E}lliptic
  {P}roblems.
\newblock SIAM J. Sci. Comput. \textbf{35}(4), A2025--A2045 (2013)

\bibitem{Loeper2005}
Loeper, G., Rapetti, F.: Numerical solution of the {M}onge-{A}mp\`ere equation
  by a {N}ewton's algorithm.
\newblock C. R. Math. Acad. Sci. Paris \textbf{340}(4), 319--324 (2005)

\bibitem{Mohammadi2007}
Mohammadi, B.: Optimal transport, shape optimization and global minimization.
\newblock C. R. Math. Acad. Sci. Paris \textbf{344}(9), 591--596 (2007)

\bibitem{Neilan2010}
Neilan, M.: A nonconforming {M}orley finite element method for the fully
  nonlinear {M}onge-{A}mp\`ere equation.
\newblock Numer. Math. \textbf{115}(3), 371--394 (2010)

\bibitem{Oberman2008}
Oberman, A.M.: Wide stencil finite difference schemes for the elliptic
  {M}onge-{A}mp\`ere equation and functions of the eigenvalues of the
  {H}essian.
\newblock Discrete Contin. Dyn. Syst. Ser. B \textbf{10}(1), 221--238 (2008)

\bibitem{Oliker1988}
Oliker, V.I., Prussner, L.D.: On the numerical solution of the equation
  {$(\partial\sp 2z/\partial x\sp 2)(\partial\sp 2z/\partial y\sp
  2)-((\partial\sp 2z/\partial x\partial y))\sp 2=f$} and its discretizations.
  {I}.
\newblock Numer. Math. \textbf{54}(3), 271--293 (1988)

\bibitem{Ong94}
Ong, M.E.G.: Uniform refinement of a tetrahedron.
\newblock SIAM J. Sci. Comput. \textbf{15}(5), 1134--1144 (1994)

\bibitem{Vuong10}
Vuong, A.V., Heinrich, C., Simeon, B.: I{SOGAT}: a 2{D} tutorial {MATLAB} code
  for isogeometric analysis.
\newblock Comput. Aided Geom. Design \textbf{27}(8), 644--655 (2010)

\bibitem{Wang2004}
Wang, X.J.: On the design of a reflector antenna. {II}.
\newblock Calc. Var. Partial Differential Equations \textbf{20}(3), 329--341
  (2004)

\bibitem{Zheligovsky10}
Zheligovsky, V., Podvigina, O., Frisch, U.: The {M}onge-{A}mp\`ere equation:
  various forms and numerical solution.
\newblock J. Comput. Phys. \textbf{229}(13), 5043--5061 (2010)

\end{thebibliography}
\end{document}